\documentclass[11pt]{amsart}
\usepackage{amsmath,amssymb,mathtools,booktabs}
\usepackage[colorlinks=true,linkcolor=blue,urlcolor=blue,citecolor=blue]{hyperref}
\newtheorem{theorem}{Theorem}[section]
\newtheorem{lemma}[theorem]{Lemma}
\newtheorem{proposition}[theorem]{Proposition}

\theoremstyle{remark}
\newtheorem{remark}[theorem]{Remark}

\newcommand{\D}{\mathbb{D}}
\newcommand{\C}{\mathbb{C}}
\newcommand{\T}{\mathbb{T}}
\newcommand{\R}{\mathbb{R}}
\newcommand{\Z}{\mathbb{Z}}
\newcommand{\A}{\mathcal{A}}
\newcommand{\Hs}{\mathcal{H}}
\newcommand{\cstar}{0.52755}
\DeclareMathOperator{\artanh}{artanh}
\DeclareMathOperator{\Real}{Re}
\DeclareMathOperator{\Imag}{Im}
\DeclarePairedDelimiter{\abs}{\lvert}{\rvert}
\DeclarePairedDelimiter{\norm}{\lVert}{\rVert}

\title[Contractive multipliers and Korenblum's constant]
{Contractive multipliers on annular Hardy spaces and Korenblum's constant}

\author{Frank Wikstr\"om}
\address{Centre for Mathematical Sciences, Lund University,
Box 118, SE-221 00 Lund, Sweden}
\email{frank.wikstrom@math.lth.se}

\begin{document}
\begin{abstract}
Korenblum's maximum principle asserts that there is a constant $c\in(0,1)$ such
that, if $f$ and $g$ in the Bergman space $A^2(\D)$ satisfy $\abs{f}\le\abs{g}$
on the annulus $c<\abs{z}<1$, then $\norm{f}_{A^2}\le\norm{g}_{A^2}$. We prove
that the best such constant satisfies $c_2\ge0.52755$, improving the previous
bounds, which are at most $0.447$. The main new tool is a family of coefficient
inequalities. If $\abs{f}\le\abs{g}$ on the annulus, then $f/g$ is a contractive
multiplier on every weighted Hardy space of the annulus, and each such space
gives an inequality between the Taylor coefficients of $f$ and $g$. These
inequalities alone give $c_2\ge0.52$ through an explicit Stieltjes
interpolation. Combined with a sharp form of Wang's circle estimate and a
computer-assisted certificate, they give $0.52755$. An exact rational example
shows that this combined criterion cannot reach $0.527558$.
\end{abstract}
\maketitle

\section{Introduction}
Let $\D$ be the unit disk and let $dA(z)=\pi^{-1}\,dx\,dy$ be
normalized area measure. The Bergman space $A^2(\D)$ consists of the
holomorphic functions $f$ on $\D$ for which
\[
 \norm{f}_{A^2}^2 = \int_{\D} \abs{f(z)}^2\,dA(z) < \infty.
\]
Korenblum's maximum principle asserts that there exists an absolute
constant $c\in(0,1)$ with the following property: if $f, g\in A^2(\D)$
and
\begin{equation}\label{eq:domination}
    \abs{f(z)} \leq \abs{g(z)},\qquad \text{for $c < \abs{z} < 1$},
\end{equation}
then $\norm{f}_{A^2}\leq \norm{g}_{A^2}$. We call such numbers $c$
\emph{admissible}. Every number below an admissible one is admissible, and the
supremum of the admissible numbers is denoted by $c_2$ and is called Korenblum's
constant.

The difficulty lies in the zeros of $g$ in the disk $\abs{z}\le c$. If $g$ has
no zeros there, then $f/g$ is holomorphic in $\D$, since every zero of $g$ in
the annulus is a zero of $f$ of at least the same order. The maximum principle
then gives $\abs{f}\le\abs{g}$ in all of~$\D$, and the conclusion is
immediate. On the other hand, the pair $f\equiv a$, $g(z)=z$
satisfies~\eqref{eq:domination} with $c=a$, and
$\norm{f}_{A^2}^2=a^2>1/2=\norm{g}_{A^2}^2$ when $a>1/\sqrt2$. Hence
$c_2\le1/\sqrt2$.

The maximum principle was conjectured by Korenblum~\cite{Korenblum1991} and
proved by Hayman~\cite{Hayman1999} with $c=0.04$. The lower bound was then
improved by Hinkkanen~\cite{Hinkkanen1999} to $0.15724$, by
Schuster~\cite{Schuster2006} to $0.21$, and by Wang~\cite{Wang2006,Wang2011}
to $0.25018$ and $0.28185$. Wang later obtained $c_2 \ge 0.3554$ using a
coefficient estimate on the circles of the annulus~\cite{Wang2025}. By
combining Wang's estimate with a moment-duality criterion, the
author~\cite{Wikstrom2026} improved the bound to $c_2 \ge 0.4263$.
Unrefereed computer-assisted certificates in Tao's repository of optimization
constants~\cite{TaoRepo} refine the circle estimate within this criterion and
reach $0.432$~\cite{Pinho2026} and $0.44605$~\cite{Patankar2026}, and $0.447$
with an extra correction at the first Taylor coefficient~\cite{Pinho2026b}. In
the other direction, the
best known upper bound is $c_2<0.6778995$, due to Wang~\cite{Wang2008}.

The main result of this paper is the following improvement of the lower bound.

\begin{theorem}\label{thm:main} The number $\cstar$ is admissible. Equivalently,
if $f,g \in A^2(\D)$ and $\abs{f(z)} \le \abs{g(z)}$ for $\cstar < \abs{z} < 1$, then
$\norm{f}_{A^2} \le \norm{g}_{A^2}$. In particular $c_2 \ge \cstar$.
\end{theorem}

\subsection*{The kernel inequalities}

Theorem~\ref{thm:main} rests on a new family of coefficient inequalities. If
$\abs{f}\le\abs{g}$ on $\A=\{c<\abs{z}<1\}$, then $\omega=f/g$ is holomorphic on
$\A$ and $\abs{\omega}\le1$. If $\omega$ were a contractive multiplier on
$A^2(\D)$, there would be nothing to prove, but~$\omega$ is controlled only on
$\A$, and~$\omega$ may have poles in $\abs{z}<c$. Instead, we use Hilbert spaces
of functions on $\A$ on which the hypothesis acts as a contraction. For $s>0$,
let $\Hs_s$ be the Hardy space of $\A$, normed by the mean square on the outer
boundary circle plus~$s$ times the mean square on the inner one. Then $\omega$
is a contractive multiplier on $\Hs_s$. Up to normalization, the spaces $\Hs_s$
are the Hardy spaces of character-automorphic functions on $\A$ that appear in
Abrahamse's Nevanlinna--Pick theorem for the annulus~\cite{Abrahamse1979}
(Remark~\ref{rem:positive-kernel}). Put $q=c^2$ and
\[
 h_s(k)=\frac{(1+s)q^k}{q^k+s},\qquad s>0,\ k\ge 0.
\]

\begin{theorem}\label{thm:kernel}
Let $f=\sum a_kz^k$ and $g=\sum b_kz^k$ be holomorphic in $\D$ with
$\abs{f}\le\abs{g}$ on $\A$. Then, for every $s>0$,
\begin{equation}\label{eq:kernel-coeff}
 \sum_{k\ge 0}\bigl(\abs{b_k}^2-\abs{a_k}^2\bigr) h_s(k) \ge 0,
\end{equation}
and the series converges absolutely.
\end{theorem}

The proof in Section~\ref{sec:kernels} is short. The adjoint of multiplication
by $\omega$ on $\Hs_s$ maps an explicit vector built from the coefficients of
$g$ to the corresponding vector built from $f$, and the squared norms of these
vectors are, up to the factor $1+s$, the two sums in~\eqref{eq:kernel-coeff}.
Equivalently,~\eqref{eq:kernel-coeff} is the average over a circle of a
positive kernel, which explains the name \emph{kernel inequalities}. No
assumption is made on the zeros of $g$ in $\abs{z}\le c$, and the Bergman space
enters only later, through the weights $1/(k+1)$.

As $s\to0$ and $s\to\infty$, the weights $h_s(k)$ tend to $1$ and $q^k$ for
each fixed $k$, which correspond to the comparison of $f$ and $g$ in mean
square on the outer and inner boundary circle, respectively. The weights
$h_s(k)$ are not averages of circle weights $t^k$ with $q\le t\le1$, whereas
every mean square comparison of $f$ and $g$ on a circle in $\A$ is an average of
the inequalities~\eqref{eq:kernel-coeff} (Remark~\ref{rem:kernel}). In this
sense the kernel inequalities use the whole annulus at once.

\subsection*{A lower bound from the kernel inequalities alone}

Suppose that $\nu$ is a positive measure on $(0,\infty)$ with
\[
 \int_0^\infty h_s(k)\,d\nu(s)=\frac1{k+1},\qquad k\ge 0.
\]
Integrating~\eqref{eq:kernel-coeff} against $\nu$ then gives
$\norm{g}_{A^2}\ge\norm{f}_{A^2}$ directly.
With $s=1/t$, such a measure corresponds to a Stieltjes function
$F(z)=\int_0^\infty d\sigma(t)/(t+z)$ with $F(q^{-k})=1/(k+1)$ for every
$k\ge 0$. In Section~\ref{sec:kernel-only} we construct such a function, with
$\sigma$ absolutely continuous with an explicit positive density.

\begin{theorem}\label{thm:kernel-only}
The number $0.52$ is admissible.
\end{theorem}

This is well above the numerical limit $0.4461$ of the circle inequalities
described below, and a simpler version of the construction with only two free
parameters already gives $0.45$, above all bounds obtained before this paper
(Remark~\ref{rem:kernel-only}). The proof uses no inequality other
than~\eqref{eq:kernel-coeff} and needs no estimates for individual Taylor
degrees. The interpolation handles all degrees at once, and the only
computation is a ball-arithmetic check that an explicit density on
$(0,\infty)$ is positive.

\subsection*{Adding the circle estimate}

Theorem~\ref{thm:main} shows how far the method can be pushed with a larger
computation. For it, we return to the moment-duality approach
of~\cite{Wikstrom2026}. Split the sequence $\abs{a_k}^2-\abs{b_k}^2$ into its
positive and negative parts $x_k$ and $y_k$. Since
$\norm{f}_{A^2}^2-\norm{g}_{A^2}^2 = \sum_k(x_k-y_k)/(k+1)$, the goal is to
show that
\[
  \sum_k \frac{x_k}{k+1}\le\sum_k\frac{y_k}{k+1}.
\]
The pointwise domination~\eqref{eq:domination} yields linear inequalities of
the form
\begin{equation}\label{eq:generic}
 \sum_{k\ge 0}x_k\,w(k)\le\sum_{k\ge 0}y_k\,w'(k),
\end{equation}
with $0\le w'\le w$. If a positive combination of such inequalities satisfies
$w(k)\ge1/(k+1)\ge w'(k)$ for every $k$, then
\[
 \sum_k\frac{x_k}{k+1}\le\sum_kx_kw(k)\le\sum_ky_kw'(k)
 \le\sum_k\frac{y_k}{k+1},
\]
which is the desired norm comparison. The kernel inequalities have this form,
with $w=w'=h_s$. Wang's estimate on the circle $\abs{z}=\sqrt{t}$ is an
inequality like~\eqref{eq:generic} with $w(k)=t^k$ and $w'(k)=\beta(t)t^k$ for
an explicit transcendental weight $\beta$ with $0\le\beta(t)<1$. The bounds
$0.4263$, $0.432$ and $0.44605$ of~\cite{Wikstrom2026,Pinho2026,Patankar2026}
use only circle inequalities of this kind. Numerically, such inequalities
cannot reach above $0.4461$, even with the sharp circle factor of
Section~\ref{sec:circle}.

We combine the two families in the moment criterion of
Theorem~\ref{thm:hybrid}, which involves a measure~$\mu$ on radii and a
measure~$\nu$ on the parameters~$s$. The circle estimate does not seem to
change where this criterion stops (see below). Its role is rather to make the
problem much easier to certify. Without it, $1/(k+1)$ must be represented
exactly by kernel weights, and the construction of
Section~\ref{sec:kernel-only} becomes harder as~$c$ grows. A factor $\beta<1$
opens a gap between the two sides of the criterion, so that measures with
finitely many atoms can work. For Theorem~\ref{thm:main}, ball arithmetic
verifies the necessary inequalities for such measures up to a certain degree,
and analytic tail estimates cover the remaining degrees
(Section~\ref{sec:certificate}).

For the circle estimate, we use the sharp per-circle factor
\[
 G'_c(\rho)=\frac12\left(1+\frac1{\sqrt{1-K_c(\rho)^2}}\right),
\]
where $K_c(\rho)=c^*_{\A}(\rho,-\rho)$, and $c^*_\A$ is the M\"obius
pseudodistance of $\A$, in place of Wang's factor
\[
 G_c(\rho)=\frac12\left(1+
 \sqrt{\frac{1+3K_c(\rho)^2}{1-K_c(\rho)^2}}\right)
\]
from \cite[(8)--(9)]{Wang2025}. Then $G'_c\le G_c$, and
$G_c/G'_c\to2$ as $K_c\to1$. Section~\ref{sec:circle} gives a short proof
of the estimate with $G'_c$, based on Schuster's formula for $c^*_\A$, a
theta-function identity and a rank-two determinant bound, and shows that
$G'_c$ cannot be improved (Proposition~\ref{prop:sharp}). In the circle-only
criterion of \cite{Wikstrom2026}, replacing $G_c$ by $G'_c$ raises the
numerical limit from $0.4264$ to $0.4461$. In the present setup, the gain is
only about $10^{-5}$ (Section~\ref{sec:certificate}).

\subsection*{The limit of the method}

Theorem~\ref{thm:main} is close to the limit of what can be achieved with the
present criterion. Proposition~\ref{prop:obstruction} produces $21$ rational
numbers that satisfy every inequality used at $c=0.527558$ (and at every
larger $c<1$), while violating the Bergman norm comparison. They are
constructed from the kernel inequalities alone, and the circle inequalities
are then verified by an exact computation in the Bernstein basis. In
particular, the kernel inequalities alone give nothing beyond $0.527558$, and
numerically they fail from about $0.5275572$. The numbers are not claimed to
come from a pair satisfying~\eqref{eq:domination}, so this is not an upper
bound for $c_2$.

The obstruction also shows where information is lost. The hypothesis
$\abs{f}\le\abs{g}$ on $\A$ is equivalent to the positivity of the kernels
behind Theorem~\ref{thm:kernel}, and every inequality for the numbers
$\abs{a_k}^2$ and $\abs{b_k}^2$ that is obtained by integrating these kernels is
a positive combination of the kernel inequalities
(Remark~\ref{rem:positive-kernel}). Any further improvement must therefore use
inequalities that are not consequences of the kernel inequalities and the circle
inequalities of the form~\eqref{eq:generic} with $w(k)=t^k$ and
$w'(k)=\beta(t)t^k$. The correction at the first Taylor coefficient
in~\cite{Pinho2026b} is of this kind. Other examples would be nonlinear
consequences of the hypothesis that combine several circles, or inequalities
that keep the products $a_j\overline{a_l}$ and $b_j\overline{b_l}$ with $j\ne l$
(Remark~\ref{rem:obstruction}).

\subsection*{Organization}
Section~\ref{sec:kernels} proves Theorem~\ref{thm:kernel}, and
Section~\ref{sec:kernel-only} shows that the kernel inequalities alone give
Theorem~\ref{thm:kernel-only}. Section~\ref{sec:circle} proves the sharp
circle estimate, and Section~\ref{sec:criterion} combines the two families
into the moment criterion. Section~\ref{sec:certificate} gives the
certificate for Theorem~\ref{thm:main} and its verification.
Section~\ref{sec:obstruction} shows that the criterion fails for
$c\ge 0.527558$.

\subsection*{Notation}

Throughout the paper, $0 < c < 1$, $q=c^2$, $\A = \{ z : c< \abs{z} <1 \}$ and
$\T_\rho=\{z: \abs{z} =\rho\}$. For holomorphic functions $f(z)=\sum_{k\ge
0}a_kz^k$ and $g(z)=\sum_{k\ge 0}b_kz^k$ on $\D$, we write
\begin{equation}\label{eq:coefficient-split}
  \begin{aligned}
    x_k&=(\abs{a_k}^2-\abs{b_k}^2)_+,&
    X(t)&=\sum_{k\ge 0}x_kt^k, \\
    y_k&=(\abs{b_k}^2-\abs{a_k}^2)_+,&
    Y(t)&=\sum_{k\ge 0}y_kt^k,
 \end{aligned}
\end{equation}
where $v_+=\max\{v,0\}$. Both sums include the constant term and converge for $0
\le t < 1$. For $f,g\in A^2(\D)$,
\begin{equation}\label{eq:bergman-difference}
  \norm{f}_{A^2}^2 - \norm{g}_{A^2}^2
  =\sum_{k\ge 0} \frac{x_k-y_k}{k+1}.
\end{equation}
Circle means are
$M_h(\rho)^2= \frac{1}{2\pi} \int_0^{2\pi}\abs{h(\rho e^{i\theta})}^2\,d\theta$, and
all integrals over circles below use normalized angular measure.

\section{Coefficient inequalities from annular Hardy kernels}
\label{sec:kernels}

Let $\Hs_s$ be the space of Laurent series $F(z)=\sum_{n\in\Z} F_n z^n$ with
\[
 \norm{F}_{\Hs_s}^2=\sum_{n\in\Z} (1+sq^n)\abs{F_n}^2<\infty.
\]
Such a series converges on the annulus $\A = \{ z : c< \abs{z} <1 \}$, and
\[
 \norm{F}_{\Hs_s}^2 = \lim_{r\to 1} M_F(r)^2 + s\lim_{r\to c} M_F(r)^2,
\]
so $\Hs_s$ is the Hardy space of $\A$ with the weights $1$ and $s$ on the two
boundary circles; compare \cite{Sarason1965}. It has the reproducing kernel
\begin{equation}\label{eq:hardy-kernel}
  k_s(z,w) = \sum_{n\in\Z} \frac{(z\overline w)^n}{1+sq^n},
\end{equation}
which converges absolutely for $z, w\in\A$: the terms with $n \ge 0$ are bounded
by $\abs{z\overline w}^n$, and those with $n < 0$ by $s^{-1}(q/\abs{z\overline
w})^{\abs{n}}$.

We now prove Theorem~\ref{thm:kernel}. Recall that
$h_s(k)=(1+s)q^k/(q^k+s)$. Thus, $h_s(0)=1$, $h_s(k) \to 1$ as $s \to 0$, and
$h_s(k) \to q^k$ as $s \to \infty$. In the
notation~\eqref{eq:coefficient-split}, inequality~\eqref{eq:kernel-coeff} says
that
\[
 \sum_{k\ge0} x_k\,h_s(k) \le \sum_{k\ge0} y_k\,h_s(k).
\]

\begin{proof}[Proof of Theorem~\ref{thm:kernel}]
The series converges absolutely, since $0<h_s(k)\le(1+s)q^k/s$ and
$\sum_k(\abs{a_k}^2+\abs{b_k}^2)q^k = M_f(c)^2+M_g(c)^2$. We may assume
$g\not\equiv 0$. If $g$ has a zero in $\A$, the domination forces $f$ to vanish there to at
least the same order. Hence, $\omega=f/g$ is holomorphic on $\A$ and
$\abs{\omega}\le1$.

We claim that $\omega$ is a contractive multiplier on $\Hs_s$. Let $F$ be
a Laurent polynomial, and write $(\omega F)_n$ for the Laurent coefficients
of $\omega F$ on $\A$. For $c<r<1$ we have $\abs{\omega F} \le \abs{F}$ on
$\T_r$, and Parseval's identity gives
$\sum_{n} \abs{(\omega F)_n}^2r^{2n}\le\sum_{n} \abs{F_n}^2r^{2n}$. Letting
$r\to1$ and $r\to c$ and using Fatou's lemma, we obtain
$\sum_{n} \abs{(\omega F)_n}^2\le\sum_{n}\abs{F_n}^2$ and
$\sum_{n} \abs{(\omega F)_n}^2 q^n \le \sum_{n} \abs{F_n}^2 q^n$. Adding
these with the weights~$1$ and~$s$ gives $\norm{\omega F}_{\Hs_s} \le
\norm{F}_{\Hs_s}$. Laurent polynomials are dense in $\Hs_s$, so $F\mapsto\omega
F$ extends to a contraction~$M_\omega$ on~$\Hs_s$, and the continuity of point
evaluations shows that $M_\omega F = \omega F$ for every $F\in\Hs_s$.

To obtain the coefficient inequality, we apply $M_\omega^*$ to an explicit
vector. Put
\[
 g^\sharp(z)=\sum_{j\ge 0}\frac{\overline{b_j}}{1+sq^{-j}}\,z^{-j},
 \qquad
 f^\sharp(z)=\sum_{j\ge 0}\frac{\overline{a_j}}{1+sq^{-j}}\,z^{-j}.
\]
Since $1/(1+sq^{-j})=q^j/(q^j+s)$,
\[
 \norm{g^\sharp}_{\Hs_s}^2=\sum_{j\ge 0}\abs{b_j}^2\frac{q^j}{q^j+s}
 \qquad\text{and}\qquad
 \norm{f^\sharp}_{\Hs_s}^2=\sum_{j\ge 0}\abs{a_j}^2\frac{q^j}{q^j+s},
\]
and both are finite by the first paragraph of the proof. Let $F$ be a Laurent
polynomial and fix $c<r<1$. The Laurent series of $g$ and of $\omega F$
converge absolutely on $\T_r$, and $g\omega=f$ on $\A$, so
\[
 \langle M_\omega F,g^\sharp\rangle_{\Hs_s}
 =\sum_{j\ge 0}b_j(\omega F)_{-j}
 =\int_{\T_r}g\,\omega F
 =\int_{\T_r}fF
 =\sum_{j\ge 0}a_jF_{-j}
 =\langle F,f^\sharp\rangle_{\Hs_s}.
\]
Since Laurent polynomials are dense, $M_\omega^*g^\sharp=f^\sharp$. Hence
$\norm{f^\sharp}_{\Hs_s}\le\norm{g^\sharp}_{\Hs_s}$, that is,
\[
 \sum_{k\ge 0}\bigl(\abs{b_k}^2-\abs{a_k}^2\bigr)\frac{q^k}{q^k+s} \ge 0.
\]
Multiplication by $1+s$ gives~\eqref{eq:kernel-coeff}.
\end{proof}

\begin{remark}\label{rem:positive-kernel}
(a) The proof does not depend on the radius $r$, and no assumption is made on
the zeros of $g$ in $\abs{z}\le c$.

(b) The inequalities~\eqref{eq:kernel-coeff} can also be read off from a
positive kernel, which explains their name. Since
$M_\omega^*k_s(\cdot,w)=\overline{\omega(w)}\,k_s(\cdot,w)$, the kernel
$\bigl(1-\omega(z)\overline{\omega(w)}\bigr)k_s(z,w)$ is positive semidefinite
on~$\A$ (see, e.g., \cite{AglerMcCarthy2002}), and hence so is
\[
 P_s(z,w)=\bigl(g(z)\overline{g(w)}-f(z)\overline{f(w)}\bigr)k_s(z,w).
\]
Conversely, $P_s(z,z)\ge0$ gives $\abs{f(z)}\le\abs{g(z)}$, so for each $s$
the positivity of $P_s$ is equivalent to $\abs{f}\le\abs{g}$ on~$\A$. For every
$c<\rho<1$ we have $g^\sharp=\int_{\T_\rho}\overline{g(w)}\,k_s(\cdot,w)$, and
similarly for $f^\sharp$. So $\norm{g^\sharp}^2-\norm{f^\sharp}^2$ is the
average of $P_s$ over $\T_\rho\times\T_\rho$.

(c) For $0\le\alpha<1$ and $s=q^\alpha$, the map $F\mapsto z^\alpha F$ is an
isometry from $\Hs_s$ onto the Hardy space of multivalued holomorphic functions
on $\A$ that are multiplied by $e^{2\pi i\alpha}$ under continuation once
around $\A$, normed by the sum of the mean squares on the two boundary
circles. Thus, up to factors $\varphi(z)\overline{\varphi(w)}$, the kernels
$k_s$ with $q<s\le1$ are the reproducing kernels of the Hardy spaces of
character-automorphic functions on $\A$ (compare \cite{Sarason1965}). These
are the kernels in Abrahamse's Nevanlinna--Pick theorem for finitely connected
domains~\cite{Abrahamse1979}; see also~\cite{VinnikovFedorov1998}. Abrahamse
normalizes the boundary measure differently, but changing the boundary
measure by a smooth positive density only shifts the character and multiplies
the kernels by factors $\varphi(z)\overline{\varphi(w)}$. The other values of
$s$ give nothing new in this respect, since
$k_{sq^m}(z,w)=(z\overline w)^{-m}k_s(z,w)$ for $m\in\Z$.

(d) More generally, let $\lambda$ be a complex measure with compact support
in~$\A$, and put $v_g=\int\overline{g(w)}\,k_s(\cdot,w)\,d\lambda(w)$ and
$v_f=\int\overline{f(w)}\,k_s(\cdot,w)\,d\lambda(w)$. Then
$M_\omega^*v_g=v_f$, so $\norm{v_f}\le\norm{v_g}$. This is the integral of
$P_s$ against $d\lambda(w)\,\overline{d\lambda(z)}$, and it involves the
products $b_j\overline{b_l}$ and $a_j\overline{a_l}$ with $j\ne l$. Averaging
over the rotated measures $\lambda(e^{-it}\,\cdot\,)$, $0\le t<2\pi$, removes
these products, and a short computation gives
\[
 \sum_{k\ge 0}\bigl(\abs{b_k}^2-\abs{a_k}^2\bigr)
 \sum_{m\in\Z}\abs{\lambda_m}^2\,\frac{q^k}{q^k+sq^m}\ge0,
 \qquad\text{where }
 \lambda_m=\int\overline{w}^{\,m}\,d\lambda(w).
\]
This is a positive combination of the inequalities~\eqref{eq:kernel-coeff}
with $s$ replaced by $sq^m$. Normalized arc length on $\T_\rho$ gives
Theorem~\ref{thm:kernel}, and $w^m$ times it gives the parameter $sq^m$. In
particular, every inequality that involves only the numbers $\abs{a_k}^2$ and
$\abs{b_k}^2$ and is obtained by integrating $P_s$ against sums of such
measures $d\lambda(w)\,\overline{d\lambda(z)}$ is a positive combination of the
inequalities~\eqref{eq:kernel-coeff}.
\end{remark}

The following representations are also used in Section~\ref{sec:obstruction}.

\begin{lemma}\label{lem:stieltjes}
Let $0<\alpha<1$, $s>0$ and $x>0$. Then
\begin{align}
 x^\alpha&=\frac{\sin\pi\alpha}{\pi}\int_0^\infty
 \frac{x}{x+y}\,y^{\alpha-1}\,dy,\label{eq:stieltjes-power}\\
 \frac{x^\alpha}{x^\alpha+s}&=\frac{s\sin\pi\alpha}{\pi}\int_0^\infty
 \frac{x}{x+y}\,
 \frac{y^{\alpha-1}}{y^{2\alpha}+2sy^\alpha\cos\pi\alpha+s^2}\,dy.
 \label{eq:stieltjes-ratio}
\end{align}
Both densities are positive.
\end{lemma}

\begin{proof}
Identity~\eqref{eq:stieltjes-power} is the beta integral $\int_0^\infty
y^{\alpha-1}(x+y)^{-1}\,dy=\pi x^{\alpha-1}/\sin\pi\alpha$. For
\eqref{eq:stieltjes-ratio}, apply Cauchy's integral formula to
$z^{\alpha-1}/(z^\alpha+s)$ on the slit plane $\C\setminus(-\infty,0]$. The
denominator does not vanish there, since $\abs{\arg z^\alpha} < \pi\alpha < \pi$. The
density is $-\pi^{-1}$ times the imaginary part of the boundary value at
$-y+i0$, which is the stated expression. It is $O(y^{\alpha-1})$ at $0$ and
$O(y^{-\alpha-1})$ at $\infty$, so the integral converges; the contributions of
small and large circles vanish for the same reason.
\end{proof}

\begin{remark}\label{rem:kernel}
For $q<t<1$ put $\alpha=\log t/\log q\in(0,1)$. With $x=q^k$,
\eqref{eq:stieltjes-power} gives
\[
 t^k=\frac{\sin\pi\alpha}{\pi}\int_0^\infty h_y(k)\,
 \frac{y^{\alpha-1}}{1+y}\,dy.
\]
So the mean square comparison $M_f(\sqrt t)\le M_g(\sqrt t)$, that is
$\sum_k(\abs{b_k}^2-\abs{a_k}^2) t^k \ge 0$, is an average of the
inequalities~\eqref{eq:kernel-coeff}. The converse fails. If
$h_s(k)=\int_{[q,1]}t^k\,d\mu(t)$ for a positive measure $\mu$, then
$\mu([q,1])=h_s(0)=1$. Moreover, $h_s(k)/q^k=\int(t/q)^k\,d\mu$ tends to the
finite limit $(1+s)/s$, so $\mu$ is concentrated at $q$. But then
$h_s(k)/q^k\le1$, whereas $h_s(k)/q^k=(1+s)/(q^k+s)>1$ for $k\ge1$.
\end{remark}

\section{A lower bound from the kernel inequalities alone}\label{sec:kernel-only}
This section proves Theorem~\ref{thm:kernel-only}. The first step is the
observation from the introduction.

\begin{proposition}\label{prop:average}
Let $\nu$ be a positive measure on $(0,\infty)$ such that
\begin{equation}\label{eq:interpolation}
 \int h_s(k)\,d\nu(s)=\frac1{k+1}, \qquad (k=0,1,2,\dots).
\end{equation}
If $f,g\in A^2(\D)$ satisfy $\abs{f(z)} \le \abs{g(z)}$ for all $c < \abs{z} <
1$, then
\[
 \norm{g}_{A^2}^2-\norm{f}_{A^2}^2
 =\int\sum_{k\ge 0}\bigl(\abs{b_k}^2-\abs{a_k}^2\bigr)h_s(k)\,d\nu(s) \ge 0.
\]
In particular, $c$ is admissible.
\end{proposition}

Note that condition~\eqref{eq:interpolation} depends on $q = c^2$, even if
the notation does not show it explicitly.

\begin{proof}
Put $d_k=\abs{b_k}^2-\abs{a_k}^2$. By~\eqref{eq:interpolation} and
Tonelli's theorem,
\[
 \int\sum_{k\ge 0}\abs{d_k}h_s(k)\,d\nu(s)=
 \sum_{k\ge 0}\frac{\abs{d_k}}{k+1}
 \le \norm{f}_{A^2}^2+\norm{g}_{A^2}^2.
\]
Hence Fubini's theorem allows us to interchange sum and integral, and
\eqref{eq:interpolation} gives the identity. The integrand is non-negative
by Theorem~\ref{thm:kernel}.
\end{proof}

Taking $k=0$ in~\eqref{eq:interpolation} shows that $\nu$ is a probability
measure.

\subsection{Stieltjes functions}

Let $\Omega = \C \setminus (-\infty,0]$, with the principal branches of $\sqrt
z$ and $\log z$. For $t>0$, we have $h_{1/t}(k)=(1+t)/(t+q^{-k})$. Hence, if
$\rho\ge 0$ on $(0,\infty)$ and $\nu$ is the image of $\rho(t)\,dt/(1+t)$ under
$t\mapsto1/t$, then
\[
 \int h_s(k)\,d\nu(s)=\int_0^\infty\frac{\rho(t)}{t+q^{-k}}\,dt.
\]
So~\eqref{eq:interpolation} holds as soon as the Stieltjes function
$F(z)=\int_0^\infty\rho(t)\,dt/(t+z)$ satisfies $F(q^{-k})=1/(k+1)$ for
every $k\ge 0$. The density is recovered from the boundary values of $F$ on
the cut.

\begin{lemma}\label{lem:stieltjes-rep}
Let $F$ be holomorphic on $\Omega$ with $F(\bar z)=\overline{F(z)}$, bounded
on $\Omega\cap\{\abs{z} < 1\}$, and with $F(z)\to0$ as $z\to\infty$ in
$\Omega$. Suppose that the restriction of $F$ to the upper half-plane extends
continuously to $(-\infty,0)$, with boundary values $F(-t+i0)$, and that
\[
 \rho(t)=-\frac1\pi\Imag{F}(-t+i0)\ge 0\qquad\text{(for $t>0$)}.
\]
Then $\int_0^\infty\rho(t)\,dt/(1+t)<\infty$, and
$F(z)=\int_0^\infty\rho(t)\,dt/(t+z)$ for every $z\in\Omega$.
\end{lemma}

\begin{proof}
Fix $z\in\Omega$ and $0<\varepsilon<\abs{z}<R$. Since the boundary values are
continuous, Cauchy's formula holds on $\Omega\cap\{\varepsilon<\abs\zeta<R\}$
with the two edges of $[-R,-\varepsilon]$ as part of the boundary. By symmetry,
the boundary values from below are $\overline{F(-t+i0)}$, so the two edges
contribute
\[
 \frac1{2\pi i}\int_\varepsilon^R\frac{F(-t-i0)-F(-t+i0)}{t+z}\,dt
 =\int_\varepsilon^R\frac{\rho(t)}{t+z}\,dt.
\]
The circles $\abs\zeta=R$ and $\abs\zeta=\varepsilon$ contribute at most
$R\sup_{\abs\zeta=R}\abs F/(R-\abs{z})$ and
$\varepsilon\sup_{\abs\zeta=\varepsilon}\abs F/(\abs{z}-\varepsilon)$, which
tend to $0$ as $R\to\infty$ and $\varepsilon\to 0$. Thus
$F(z)=\lim\int_\varepsilon^R\rho(t)\,dt/(t+z)$. For $z=1$, monotone convergence
gives $\int_0^\infty\rho(t)\,dt/(1+t)=F(1)<\infty$. For fixed $z\in\Omega$,
the ratio $(1+t)/\abs{t+z}$ is bounded for $t>0$, so dominated convergence gives
the formula for every $z\in\Omega$.
\end{proof}

\subsection{An explicit interpolant}

For $z\in\Omega$, put
\[
 \eta(z)=\frac{\sqrt z-\sqrt q}{\sqrt z+\sqrt q},\qquad
 \mathcal B(z)=\prod_{n\ge 0}\frac{1-q^{n/2}\sqrt z}{1+q^{n/2}\sqrt z}.
\]
Since $\Real\sqrt z>0$ on $\Omega$, we have $\abs\eta<1$. Also, every factor in
the definition of $\mathcal B$ has modulus less than~$1$, and the product
converges locally uniformly. Thus, $\abs{\mathcal B}<1$ on $\Omega$, $\mathcal B$
vanishes at the point $q^{-n}$, for each $n\ge 0$, and
\[
 \mathcal B(q)=\prod_{n\ge1}\frac{1-q^{n/2}}{1+q^{n/2}}>0.
\]
For a real polynomial $C(\eta)=\sum_{j=0}^MC_j\eta^j$ with $C_0=1$, define
\begin{equation}\label{eq:interpolant}
 F(z)=\log\frac{1}{q}\,\biggl(\frac{1}{\log(z/q)}
 -\frac{q}{\mathcal B(q)}\,\frac{\mathcal B(z)\,C(\eta(z))}{z-q}\biggr).
\end{equation}
The function $\log(1/q)/\log(z/q)$ takes the value $1/(k+1)$ at $z=q^{-k}$, but
it has a pole at $z=q$, so it is not a Stieltjes function. The second term in
\eqref{eq:interpolant} removes the pole without changing the values at the
points $q^{-k}$, where~$\mathcal B$ vanishes. We view the polynomial~$C$ as a
free parameter, to be chosen so that the density of~$F$ is positive.

\begin{lemma}\label{lem:interpolant}
The function $F$ in~\eqref{eq:interpolant} satisfies $F(q^{-k})=1/(k+1)$ for
every $k\ge 0$, and it satisfies the hypotheses of
Lemma~\ref{lem:stieltjes-rep}, except possibly $\rho\ge 0$. For $t>0$ we have
$\rho(t)=\pi^{-1}\log(1/q)\,\Psi(t)$, where
\begin{equation}\label{eq:density}
 \Psi(t)=\frac\pi{\log^2(t/q)+\pi^2}
 -\frac
  {q\Imag\bigl(\mathcal B_+(t)\,C(\eta_+(t))\bigr)}
  {(t+q)\,\mathcal B(q)}
\end{equation}
and
\[
  \mathcal B_+(t) = \prod_{n\ge 0}
  \frac{1-iq^{n/2}\sqrt t}{1+iq^{n/2}\sqrt t},
  \qquad
  \eta_+(t)=\frac{i\sqrt t-\sqrt q}{i\sqrt t+\sqrt q}.
\]
Every factor of $\mathcal B_+(t)$, and $\eta_+(t)$ itself, has modulus $1$.
\end{lemma}

\begin{proof}
In $\Omega$, $\log(z/q)$ vanishes only at $z=q$, where
$1/\log(z/q)=q/(z-q)+O(1)$. Since $\eta(q)=0$ and $C(0)=1$, also
$\mathcal B(z)C(\eta(z))/\mathcal B(q)=1+O(z-q)$. So the poles of the two
terms cancel, and $F$ is holomorphic on $\Omega$. At $z=q^{-k}$ we have
$\log(z/q)=(k+1)\log(1/q)$ and $\mathcal B(z)=0$, so $F(q^{-k})=1/(k+1)$.
All ingredients are real on $(0,\infty)$, which gives
$F(\bar z)=\overline{F(z)}$.

Let $\norm{C}_1=\sum_j\abs{C_j}$. On $\Omega$ we have
$\abs{C(\eta)}\le\norm{C}_1$, $\abs{\mathcal B}\le1$ and
$\abs{\log(z/q)}\ge\abs{\log(\abs z/q)}$. Hence, $F(z)\to0$ as $z\to\infty$, and
$F$ is bounded on $\Omega\cap\{\abs z<q/2\}$. On the closed upper half-plane,
$\Real \sqrt z\ge 0$, so $\abs{1+q^{n/2}\sqrt z}\ge1$ and $\abs{\sqrt z+\sqrt
q}\ge\sqrt q$ there. The product for $\mathcal B$ therefore converges uniformly
on compact subsets of the closed upper half-plane minus the origin, and
$\mathcal B$ and $\eta$ extend continuously to $(-\infty,0)$ with the boundary
values $\mathcal B_+(t)$ and $\eta_+(t)$ at $-t$. Since
$\log(z/q)\to\log(t/q)+i\pi$ as $z\to-t$ from above, $F$ extends continuously as
well, and it is bounded on $\Omega\cap\{\abs z<1\}$. Taking imaginary parts
in~\eqref{eq:interpolant} gives~\eqref{eq:density}. Finally,
$\abs{1-ix}=\abs{1+ix}$ and $\abs{i\sqrt t-\sqrt q}=\abs{i\sqrt t+\sqrt q}$.
\end{proof}

\subsection{Positivity of the density}

For Theorem~\ref{thm:kernel-only} we take $q=0.52^2=0.2704$, $M=96$, and the
decimal coefficients $C_j$ listed in~\cite{Supplement}, treated as exact
rational numbers. They have at most $17$ significant digits, $C_0=1$,
$\norm{C}_1<565$ and $\norm{C'}_1=\sum_jj\abs{C_j}<18082$. By
Proposition~\ref{prop:average} and Lemmas~\ref{lem:stieltjes-rep}
and~\ref{lem:interpolant}, it suffices to show that $\Psi(t)>0$ for every $t>0$.
We write $u=\log t$ and $v=u-\log q$.

\begin{lemma}\label{lem:density-tails}
Let $H=2\bigl(\norm{C}_1/(1-\sqrt q)+\norm{C'}_1/\sqrt q\bigr)$. Then
\[
 \abs{\Imag\bigl(\mathcal B_+(t)C(\eta_+(t))\bigr)}
 \le\min\{\norm{C}_1,H\sqrt t\}.
\]
Consequently, $\Psi(t)>0$ if
\begin{equation}\label{eq:density-tails}
 He^{u/2}(v^2+\pi^2)<\pi\mathcal B(q)\qquad\text{or}\qquad
 q\norm{C}_1(v^2+\pi^2)<\pi\mathcal B(q)e^u.
\end{equation}
Both $e^{u/2}(v^2+\pi^2)$ and $e^u/(v^2+\pi^2)$ are increasing functions of
$u\in\R$.
\end{lemma}

\begin{proof}
Since $\abs{\mathcal B_+}=\abs{\eta_+}=1$, we have $\abs{\mathcal
B_+C(\eta_+)}\le\norm{C}_1$. Since $C(-1)$ is real, $\Imag(\mathcal
B_+C(\eta_+)) =\Imag(\mathcal B_+C(\eta_+)-C(-1))$, and
\[
 \abs{\mathcal B_+C(\eta_+)-C(-1)}
 \le\norm{C}_1\abs{\mathcal B_+-1}+\norm{C'}_1\abs{\eta_++1},
\]
because $\abs{C'}\le\norm{C'}_1$ on the closed unit disk. For numbers $\zeta_n$
of modulus $1$, $\abs{\prod\zeta_n-1}\le\sum\abs{\zeta_n-1}$, and
$\abs{(1-ix)/(1+ix)-1}\le2x$ for $x\ge 0$. Hence, $\abs{\mathcal B_+-1}\le2\sqrt
t/(1-\sqrt q)$. Also, $\abs{\eta_++1}=2\sqrt t/\abs{i\sqrt t+\sqrt q}\le2\sqrt
t/\sqrt q$, which proves the bound. Since $q/(t+q)\le1$ and $q/(t+q)<q/t$,
either condition in~\eqref{eq:density-tails} makes the second term of $\Psi(t)$
smaller than the first. Finally, the logarithmic derivatives of
$e^{u/2}(v^2+\pi^2)$ and $e^u/(v^2+\pi^2)$ are $1/2+2v/(v^2+\pi^2)\ge1/2-1/\pi$
and $1-2v/(v^2+\pi^2)\ge1-1/\pi$.
\end{proof}

\begin{proof}[Proof of Theorem~\ref{thm:kernel-only}]
Ball arithmetic shows that the first condition in~\eqref{eq:density-tails} holds
at $u=-48$ and the second at $u=40$; the differences between the right and left
sides exceed $0.31$ and $7\cdot10^{16}$. By Lemma~\ref{lem:density-tails},
$\Psi>0$ for $u\le-48$ and for $u\ge40$.

On $[-48,40]$ we enclose $\Psi$ in ball arithmetic \cite{Johansson2017} with
$65$ decimal digits. The interval is split into $56\,587$ intervals with
dyadic endpoints, obtained by bisecting unit intervals at most $20$ times,
until the enclosure of $\Psi$ is positive. On each interval,
$\mathcal B_+(t)$ is replaced by the product of its first $100$ factors,
with error at most $2\sqrt t\,q^{50}/(1-\sqrt q)$ by the bound for products
of unimodular numbers used above. The value $C(\eta_+(t))$ is enclosed by
Taylor's formula about the midpoint of the ball for $\eta_+(t)$, with the
remainder bounded through $\norm{C''}_1=\sum_jj(j-1)\abs{C_j}$. Finally,
$\mathcal B(q)$ is enclosed by
\[
 1\ge\prod_{n>100}\frac{1-q^{n/2}}{1+q^{n/2}}
 \ge1-\frac{2q^{101/2}}{1-\sqrt q}.
\]
Every enclosure of $\Psi$ is positive, so $\rho>0$ on $(0,\infty)$. By
Lemmas~\ref{lem:stieltjes-rep} and~\ref{lem:interpolant}, the image of
$\rho(t)\,dt/(1+t)$ under $t\mapsto1/t$ is a measure $\nu$ satisfying
\eqref{eq:interpolation} for $c=0.52$. By Proposition~\ref{prop:average}, $0.52$
is admissible. The verification takes about $20$ seconds on the author's laptop, and
the program and its output are archived at~\cite{Supplement}.
\end{proof}

\begin{remark}\label{rem:kernel-only} (a) The coefficients were found by linear
programming. The function $\Psi(t)$ is affine in $(C_1,\dots,C_M)$, and we
maximized the smallest ratio of $\Psi$ to its first term over $16\,001$ equally
spaced points $u\in[-32,32]$. In floating point, the second term of $\Psi$ never
exceeds $96\%$ of the first. None of this enters the proof.

(b) In floating point, the same linear program finds positive densities up to
about $c=0.522$ with $M=96$, but none at $c=0.525$ with $M\le192$, although the
margin improves with $M$. We have not pushed the construction further. By
Remark~\ref{rem:obstruction}(a), no measure $\nu$ satisfies
\eqref{eq:interpolation} when $c\ge 0.527558$. Near that value, it is easier to
give the kernel inequalities some slack with the circle estimate, as in
Sections~\ref{sec:circle}--\ref{sec:certificate}.

(c) Only a few coefficients are needed to pass all earlier bounds. The same
verification shows that the values of $c$ in Table~\ref{tab:small-M} are
admissible, with polynomials $C$ of small degree $M$. For $M=0$ the
interpolant~\eqref{eq:interpolant} is completely explicit, and for $M=2$ we may
take $C(\eta)=1-0.97\eta+0.55\eta^2$, which gives $0.45$. The coefficients and
verification reports are archived at~\cite{Supplement}.

\begin{table}[ht]
\centering\small
\caption{Admissible values of $c$ obtained from the kernel inequalities alone,
with a polynomial $C$ of degree $M$ in~\eqref{eq:interpolant}.}
\label{tab:small-M}
\begin{tabular}{@{}lccccccc@{}}
\toprule
$M$ & $0$ & $1$ & $2$ & $4$ & $8$ & $16$ & $96$\\
\midrule
$c$ & $0.39$ & $0.41$ & $0.45$ & $0.47$ & $0.49$ & $0.50$ & $0.52$\\
\bottomrule
\end{tabular}
\end{table}
\end{remark}

\section{The sharp circle estimate}\label{sec:circle}
For $a,b\in\D$ let
\[
  \delta(a,b)=\abs[\Big]{\frac{a-b}{1-\overline b a}}
\]
be the pseudo-hyperbolic distance, and for $z,w\in\A$ let
\[
 c^*_{\A}(z,w)
 =\sup\{\delta(\psi(z),\psi(w)):\psi\in\operatorname{Hol}(\A,\D)\}
\]
be the M\"obius pseudodistance of $\A$. We will use Schuster's formula for
$c^*_\A$ on a circle. Define
\[
 E(w)=(1-w)\prod_{n\ge1}(1-q^nw)(1-\frac{q^n}{w}),\qquad \text{for $w\in\C^*$}.
\]
The product converges locally uniformly. The function $E$ has real Laurent
coefficients and satisfies the functional equations
\begin{equation}\label{eq:geo-functional}
 E(\overline w)=\overline{E(w)},\qquad
 E(1/w)=-\frac{E(w)}w = E(qw),
\end{equation}
which follow from reversing $w$ and shifting the index in the product. In
particular, $E(1)=0$, $E(-1)>0$, and
\[
 E(c)=\prod_{n\ge1}(1-c^{2n-1})^2,\qquad
 E(-c)=\prod_{n\ge1}(1+c^{2n-1})^2.
\]
Put
\begin{equation}\label{eq:geo-kappa}
 \kappa=\frac{E(c)}{E(-c)}
       =\prod_{n\ge1}\left(\frac{1-c^{2n-1}}{1+c^{2n-1}}\right)^2 \!.
\end{equation}

Schuster~\cite[pp.~3527--3528]{Schuster2006} proved the following
formula, which we state in a slightly different notation.
\begin{lemma}[Schuster's formula]\label{lem:geo-schuster}
For $c<\rho<1$ and $\abs{u} = 1$,
\begin{equation}\label{eq:geo-schuster}
  c^*_{\A}(\rho,\rho u)
  = C(\rho)\frac{\abs{E(u)}}{\abs{E(\rho^2u)}},\qquad
 \text{where } C(\rho)=\frac c\rho\,\frac{E(-\rho^2/c)}{E(-c)}.
\end{equation}
\end{lemma}

For $u=-1$,~\eqref{eq:geo-schuster} is the product formula for
$K_c(\rho)=c^*_\A(\rho,-\rho)$ used in \cite[(7)]{Wang2025} and
\cite{Wikstrom2026}. Here, we will use it for every pair of points on $\T_\rho$.

\subsection{Two theta identities}
\begin{lemma}\label{lem:theta}
For real $x\ne0$ and $\abs{u} = 1$,
\begin{align}
   E(-1)^2\abs{E(xu)}^2
    &=E(-x)^2\abs{E(u)}^2+E(x)^2\abs{E(-u)}^2,\label{eq:geo-theta1}
  \\
  E(-1)^2\,\frac qx\,E(-x/c)^2
    &=E(-c)^2E(-x)^2-E(c)^2E(x)^2.\label{eq:geo-theta2}
\end{align}
\end{lemma}

\begin{proof}
Let $\Theta_1$ and $\Theta_2$ be the spaces of holomorphic functions $F$ on
$\C^*$ with $F(qu)=F(u)/(qu^2)$ and $F(qx)=F(x)/x^2$, respectively. For a
Laurent series $F=\sum_{n} a_n u^n$, these conditions imply $a_{n+2}=q^{n+1}a_n$
and $a_{n+2}=q^na_n$ respectively, so both spaces have dimension at most two.

By~\eqref{eq:geo-functional}, the function $u\mapsto E(xu)E(x/u)$ belongs
to $\Theta_1$ for every real $x\ne0$. For $x=1$ and $x=-1$ its values at
$u=1,-1$ are $(0,E(-1)^2)$ and $(E(-1)^2,0)$, so these two functions form a
basis of $\Theta_1$. Expanding $E(xu)E(x/u)$ in this basis and evaluating at
$u=\pm1$ gives~\eqref{eq:geo-theta1}, because $E(xu)E(x/u)=\abs{E(xu)}^2$ when
$x$ is real and $\abs{u}=1$.

Similarly, $E(x)^2$ and $E(-x)^2$ belong to $\Theta_2$ and form a basis,
by their values at $x=\pm1$. The function $(q/x)E(-x/c)^2$ also belongs to
$\Theta_2$. Since $E(-1/c)=E(-c)/c$ and $E(1/c)=-E(c)/c$, its values at
$x=1$ and $x=-1$ are $E(-c)^2$ and $-E(c)^2$, which gives
\eqref{eq:geo-theta2}.
\end{proof}

Both identities are forms of the classical addition formulas for Jacobi's
theta functions.

\subsection{Closed forms on a circle}

For $c<\rho<1$ and $\theta\in\R$, put
\[
 \delta_\rho(\theta)=c^*_{\A}(\rho,\rho e^{i\theta})
 \qquad\text{and}\qquad
 S_\rho(\theta)=\frac{\delta_\rho(\theta)}
                         {\sqrt{1-\delta_\rho(\theta)^2}},
\]
so that $S_\rho(\theta)=\sinh(\frac12d_\A(\rho,\rho e^{i\theta}))$ with
$d_\A=2\artanh c^*_\A$, and $K_c(\rho)=\delta_\rho(\pi)$. For $0<x<1$ let
\[
 \chi(x)=\frac{1-x}{1+x},\qquad \Phi(x)=\prod_{n\ge 0}\chi(xq^n).
\]
Then $\kappa=\Phi(c)^2$.

\begin{lemma}\label{lem:closed}
Let $c<\rho<1$, $t=\rho^2$ and $u=e^{i\theta}$. Then
\begin{align}
 1-\delta_\rho(\theta)^2
 &=\frac{E(t)^2}{E(-1)^2}\,
   \frac{\kappa^2 \abs{E(u)}^2 + \abs{E(-u)}^2}{\abs{E(tu)}^2},\label{eq:geo-gap}\\
 \frac{S_\rho(\theta)}{S_\rho(\pi)}
 &=\frac{\kappa\abs{E(u)}}{\sqrt{\kappa^2\abs{E(u)}^2+\abs{E(-u)}^2}},
 \label{eq:S-ratio}\\
 \sqrt{1-K_c(\rho)^2}
 &=\kappa\,\frac{E(t)}{E(-t)}
 =\kappa\,\frac{1-t}{1+t}\,\Phi(qt)\,\Phi(q/t).\label{eq:stable}
\end{align}
In particular, $S_\rho(\theta)/S_\rho(\pi)$ does not depend on $\rho$.
\end{lemma}

\begin{proof}
By~\eqref{eq:geo-theta2} with $x=t$, $C(\rho)^2E(-1)^2=E(-t)^2-\kappa^2E(t)^2$.
Combining this with~\eqref{eq:geo-schuster} and~\eqref{eq:geo-theta1} for
$x=t$ gives~\eqref{eq:geo-gap}. Consequently,
\[
 S_\rho(\theta)^2=\frac{C(\rho)^2E(-1)^2}{E(t)^2}\,
 \frac{\abs{E(u)}^2}{\kappa^2\abs{E(u)}^2+\abs{E(-u)}^2},
\]
and~\eqref{eq:S-ratio} follows by putting $u=-1$ and using $E(1)=0$. Setting
$u=-1$ in~\eqref{eq:geo-gap} gives $1-K_c(\rho)^2=\kappa^2E(t)^2/E(-t)^2$. Since
$E(t)>0$ and $E(-t)>0$ for $q<t<1$, the positive square root is $\kappa
E(t)/E(-t)$, and expanding the products produces the last expression
in~\eqref{eq:stable}.
\end{proof}

All factors in the last expression of~\eqref{eq:stable} are positive, which
makes it suitable for numerical evaluation near $t=1$, where $K_c$ is close
to $1$.

\begin{lemma}[Chord inequality]\label{lem:geo-chord}
For $0\le\theta\le\pi$,
\begin{equation}\label{eq:geo-chord}
 S_\rho(\theta)\le S_\rho(\pi)\sin(\theta/2).
\end{equation}
In particular, $K_c(\rho)=\max_\theta\delta_\rho(\theta)$.
\end{lemma}

\begin{proof}
The cases $\theta=0$ and $\theta=\pi$ are clear. For $0<\theta<\pi$, since
$\abs{E(u)}=2\sin(\theta/2)\prod_n\abs{1-q^nu}^2$ and
$\abs{E(-u)}=2\cos(\theta/2)\prod_n\abs{1+q^nu}^2$,~\eqref{eq:S-ratio} shows
that~\eqref{eq:geo-chord} is equivalent to
\[
 \kappa\le
 \prod_{n\ge1}
 \frac{1+2q^n\cos\theta+q^{2n}}
      {1-2q^n\cos\theta+q^{2n}}.
\]
Each factor on the right increases with $\cos\theta$, so the right-hand side
is at least its value $\prod_{n\ge1}((1-c^{2n})/(1+c^{2n}))^2$ at
$\theta=\pi$. This exceeds~\eqref{eq:geo-kappa}, since~$\chi$ is decreasing
and $c^{2n-1}>c^{2n}$.
\end{proof}

\subsection{The estimate}
The first step holds for any pair of functions on any circle.

\begin{lemma}\label{lem:rank-two}
Let $f=\sum a_kz^k$ and $g=\sum b_kz^k$ be holomorphic in $\D$, let
$0<\rho<1$, and put
\[
 U=M_f(\rho)^2,\quad V=M_g(\rho)^2,\quad
 W=\int_{\T_\rho}f\overline g,\quad D=V-U.
\]
Then
\begin{equation}\label{eq:rank-two}
 Y(\rho^2)\le\frac{D+\sqrt{D^2+4(UV-\abs{W}^2)}}2.
\end{equation}
\end{lemma}

\begin{proof}
Let $\boldsymbol a=(a_k\rho^k)_{k\ge 0}$ and $\boldsymbol b=(b_k\rho^k)_{k\ge
0}$ in $\ell^2$, and let $T=\boldsymbol b\boldsymbol b^*-\boldsymbol
a\boldsymbol a^*$. By Parseval's identity, $\norm{\boldsymbol a}^2=U$,
$\norm{\boldsymbol b}^2=V$ and $\abs{\langle\boldsymbol a,\boldsymbol b\rangle}
= \abs{W}$. The operator $T$ has rank at most two, and since $T\le\boldsymbol
b\boldsymbol b^*$, it has at most one positive eigenvalue. Its nonzero
eigenvalues are those of the $2\times2$ matrix
\[
 \begin{pmatrix}1&0\\0&-1\end{pmatrix}
 \begin{pmatrix}V&\boldsymbol b^*\boldsymbol a\\
 \boldsymbol a^*\boldsymbol b&U\end{pmatrix},
\]
which has trace $D$ and determinant $-(UV-\abs{W}^2)$. Hence the largest
eigenvalue of $T$ is
\[
\lambda_+ = \frac12\bigl(D + \sqrt{D^2+4(UV-\abs{W}^2)} \bigr) \ge 0.
\]
Let $T_+$ be the positive part of $T$. The diagonal entries of $T$
are $(\abs{b_k}^2-\abs{a_k}^2)\rho^{2k}$, and since $T\le T_+$ and the diagonal entries
of $T_+$ are non-negative,
\[
 Y(\rho^2)=\sum_{T_{kk}>0}T_{kk}\le\sum_{T_{kk}>0}(T_+)_{kk}
 \le\operatorname{tr}T_+=\lambda_+.\qedhere
\]
\end{proof}

\begin{proposition}\label{prop:circle}
Let $f=\sum a_kz^k$ and $g=\sum b_kz^k$ be holomorphic in $\D$ with
$\abs{f}\le\abs{g}$ on $\A$. Then, for $c<\rho<1$,
\begin{equation}\label{eq:geo-estimate}
  \sum_{k\ge 0}(\abs{b_k}^2-\abs{a_k}^2)_+\rho^{2k}
  \le
  G'_c(\rho) \bigl(M_g(\rho)^2-M_f(\rho)^2 \bigr).
\end{equation}
\end{proposition}

\begin{proof}
If $g\equiv0$, or if $f$ is a unimodular constant multiple of $g$, both sides
vanish. Otherwise, $\omega=f/g$ is holomorphic on $\A$ with values in $\D$.
Let $U,V,W,D$ be as in Lemma~\ref{lem:rank-two}, and put
$L=\abs{g}^2-\abs{f}^2\ge 0$ on $\T_\rho$, so that $\int_{\T_\rho}L=D$. For
$z,w\in\T_\rho$ let $\Delta(z,w)=g(z)f(w)-f(z)g(w)$. Expanding
$\abs{\Delta}^2$ gives
\begin{equation}\label{eq:gram}
 UV-\abs{W}^2=\frac12\int_{\T_\rho}\int_{\T_\rho}\abs{\Delta(z,w)}^2.
\end{equation}
For $\zeta_1,\zeta_2\in\D$,
$1-\delta(\zeta_1,\zeta_2)^2=(1-\abs{\zeta_1}^2)(1-\abs{\zeta_2}^2)
/\abs{1-\overline{\zeta_1}\zeta_2}^2$, and therefore
\[
 \abs{\Delta(z,w)}^2 = \abs{g(z)g(w)}^2\,\abs{\omega(z)-\omega(w)}^2
 =L(z)L(w)\,\frac{\delta(\omega(z),\omega(w))^2}
                  {1-\delta(\omega(z),\omega(w))^2}.
\]
At the zeros of $g$ both sides vanish. Write $z=\rho e^{i\theta}$ and
$w=\rho e^{i\varphi}$. By the definition of $c^*_\A$ and the rotation
invariance of $\A$, $\delta(\omega(z),\omega(w))\le\delta_\rho(\theta-\varphi)$.
Lemma~\ref{lem:geo-chord}, extended to all $\theta$ by evenness and
$2\pi$-periodicity, then gives
\[
 \abs{\Delta(z,w)}^2 \le
  L(\theta)L(\varphi)\,S_\rho(\pi)^2 \sin^2\frac{\theta-\varphi}2.
\]
Since $\sin^2(x/2)=(1-\cos x)/2$ and
\[ \iint L(\theta)L(\varphi)\cos(\theta-\varphi) = \abs[\Big]{\int
L(\theta)e^{i\theta}}^2,
\]
equation~\eqref{eq:gram} yields
\[
 UV-\abs{W}^2\le\frac{S_\rho(\pi)^2}{4}
 \left(D^2-\abs[\Big]{ \int L(\theta)e^{i\theta}}^2\right)
 \le\frac{S_\rho(\pi)^2D^2}{4}.
\]
Inserting this in~\eqref{eq:rank-two} and using
$1+S_\rho(\pi)^2=(1-K_c(\rho)^2)^{-1}$ gives~\eqref{eq:geo-estimate}.
\end{proof}

The factor $G'_c(\rho)$ cannot be improved.

\begin{proposition}\label{prop:sharp}
Let $c<\rho<1$ and $G<G'_c(\rho)$. There are polynomials $f=\sum a_kz^k$ and
$g=\sum b_kz^k$ with $\abs{f}\le\abs{g}$ on $\A$ and
\[
 \sum_{k\ge 0}(\abs{b_k}^2-\abs{a_k}^2)_+\rho^{2k}
 > G\bigl(M_g(\rho)^2-M_f(\rho)^2\bigr)>0.
\]
\end{proposition}

\begin{proof}
Let $K=K_c(\rho)$, and let $0<r<1$ satisfy $2r/(1+r^2)=K$, that is,
$\delta(r,-r)=K$. Then $\sqrt{1-K^2}=(1-r^2)/(1+r^2)$, and hence
$G'_c(\rho)=1/(1-r^2)$.

By a normal family argument, there is a function
$\omega_0\in\operatorname{Hol}(\A,\D)$ with
$\delta(\omega_0(\rho),\omega_0(-\rho))=K$, and after composing with an
automorphism of $\D$ we may assume that $\omega_0(\pm\rho)=\pm r$. The odd part
$\omega(z)=\frac12(\omega_0(z)-\omega_0(-z))$ is holomorphic on $\A$, with
$\abs{\omega}<1$ and $\omega(\pm\rho)=\pm r$, and its Laurent series has only
odd powers.

Let $p_n$ be the $n$-th Fej\'er mean of the Laurent series of $\omega$. On each
circle $\T_t$, $c<t<1$, it is the convolution of $\omega$ with the Fej\'er
kernel, so $\abs{p_n}\le1$ on $\A$. Moreover, $p_n$ is an odd Laurent polynomial,
and $p_n(\pm\rho)\to\pm r$ as $n\to\infty$. Choose an integer $m$ with $2m\ge n$
and put
\[
 g_N(z)=z^{2m}\bigl(1+z^2/\rho^2\bigr)^N,\qquad f_N=p_ng_N.
\]
These are polynomials with $\abs{f_N}\le\abs{g_N}$ on $\A$. Since $g_N$ is even
and $f_N$ is odd, $b_k=0$ for odd $k$ and $a_k=0$ for even $k$. Hence, the
left-hand side equals $M_{g_N}(\rho)^2$, and
$M_{g_N}(\rho)^2-M_{f_N}(\rho)^2=(1-R_N)M_{g_N}(\rho)^2$, where
\[
 R_N=\frac{\int_0^{2\pi}\abs{p_n(\rho e^{i\theta})}^2\abs{\cos\theta}^{2N}\,d\theta}
 {\int_0^{2\pi}\abs{\cos\theta}^{2N}\,d\theta},
\]
because $\abs{g_N(\rho e^{i\theta})}=\rho^{2m}(2\abs{\cos\theta})^N$. As
$N\to\infty$, the weights concentrate at $\theta=0$ and $\theta=\pi$, and since
$\abs{p_n(-\rho)}=\abs{p_n(\rho)}$, $R_N\to\abs{p_n(\rho)}^2$. First choose $n$
with $\abs{p_n(\rho)}<1$ and $1/(1-\abs{p_n(\rho)}^2)>G$, and then $N$ with
$R_N<1$ and $1/(1-R_N)>G$.
\end{proof}

\begin{remark}\label{rem:circle}
(a) The identity
\[
 (1- \abs{\lambda}^2)T=(\boldsymbol b-\lambda\boldsymbol a)
 (\boldsymbol b-\lambda\boldsymbol a)^*
 -(\boldsymbol a-\overline\lambda\boldsymbol b)
 (\boldsymbol a-\overline\lambda\boldsymbol b)^*
\]
shows that $\lambda_+\le M_{g-\lambda f}(\rho)^2/(1- \abs{\lambda}^2)$ for every
$ \abs{\lambda}<1$, and one checks that equality holds for one $\lambda\in\D$ when
$D>0$. Wang's argument~\cite{Wang2025} bounds $Y(\rho^2)$ by
$M_{g-\lambda f}(\rho)^2/(1- \abs{\lambda}^2)$ for one particular real $\lambda$,
which leads to $G_c$.

(b) In the examples of Proposition~\ref{prop:sharp}, $L$ concentrates at
the antipodal points $\pm\rho$ and $\omega$ is extremal for
$c^*_\A(\rho,-\rho)$, so every step in the proof of
Proposition~\ref{prop:circle} is asymptotically an equality.

(c) The same factor follows from the statement that $\omega(\T_\rho)$ lies in
a hyperbolic disk of radius $\frac12d_\A(\rho,-\rho)$, combined with (a) for
$\lambda$ equal to the conjugate of the centre. That statement follows from
Lemma~\ref{lem:geo-chord} and Kirszbraun's theorem in the hyperbolic plane
\cite{Valentine1944}. Patankar~\cite{Patankar2026} also uses an enclosing
disk, with an explicit bound for its radius obtained from a flux lemma. By
Proposition~\ref{prop:sharp} his factor is at least $G'_c$, and numerically
the two agree to a relative accuracy of about $10^{-5}$ for $c\ge0.446$; in
particular, they give essentially the same limit for the circle-only
criterion. We found $G'_c$, with the argument in~(c), in September 2026,
independently of~\cite{Patankar2026}.
\end{remark}

Since $M_g(\rho)^2-M_f(\rho)^2=Y(\rho^2)-X(\rho^2)$,
Proposition~\ref{prop:circle} implies
\begin{equation}\label{eq:circle}
 X(t)\le\beta_c(t)Y(t),\qquad q\le t<1,
\end{equation}
where, by~\eqref{eq:stable},
\begin{equation}\label{eq:beta}
 \begin{aligned}
 \beta_c(t)&=1-\frac1{G'_c(\sqrt t)}=\frac{1-\sigma(t)}{1+\sigma(t)},\\
 \sigma(t)&=\sqrt{1-K_c(\sqrt t)^2}
 =\kappa\,\frac{1-t}{1+t}\,\Phi(qt)\,\Phi(q/t).
 \end{aligned}
\end{equation}
At $t=q$,~\eqref{eq:circle} holds by continuity, with $\beta_c(q)=1$. The
formula~\eqref{eq:beta} also gives $\beta_c(1)=1$, and $0\le\beta_c\le1$ on
$[q,1]$.

\section{The moment criterion}\label{sec:criterion}
For a positive measure $\nu$ on $(0,\infty)$ write
$\hat\nu(k)=\int h_s(k)\,d\nu(s)$.

\begin{theorem}\label{thm:hybrid}
Let $\mu$ be a finite positive measure on $[q,1]$ with $\mu(\{1\})=0$, and
let $\nu$ be a finite positive measure on $(0,\infty)$. If, for every integer
$k\ge 0$,
\begin{equation}\label{eq:hybrid}
 \int\beta_c(t)\,t^k\,d\mu(t)+\hat\nu(k)
 \ \le\ \frac1{k+1}\ \le\
 \int t^k\,d\mu(t)+\hat\nu(k),
\end{equation}
then $c$ is admissible.
\end{theorem}

\begin{proof}
Let $f,g\in A^2(\D)$ satisfy~\eqref{eq:domination}. Integrate
\eqref{eq:circle} against $\mu$ and~\eqref{eq:kernel-coeff} against $\nu$,
and add. Together with~\eqref{eq:hybrid} this gives
\begin{align*}
 \sum_{k\ge 0}\frac{x_k}{k+1}
 &\le\sum_{k\ge 0}x_k\Bigl(\int t^k\,d\mu+\hat\nu(k)\Bigr)\\
 &\le\sum_{k\ge 0}y_k\Bigl(\int\beta_c(t)t^k\,d\mu+\hat\nu(k)\Bigr)
 \le\sum_{k\ge 0}\frac{y_k}{k+1}.
\end{align*}
All terms are non-negative, so Tonelli's theorem justifies the interchanges,
and the last sum is at most $\norm{g}_{A^2}^2<\infty$. By
\eqref{eq:bergman-difference}, $\norm{f}_{A^2} \le \norm{g}_{A^2}$.
\end{proof}

\begin{remark}\label{rem:criterion}
(a) The proof uses only three properties of the sequences $x_k$ and $y_k$:
they are non-negative, they satisfy~\eqref{eq:circle} for $q\le t<1$, and they
satisfy~\eqref{eq:kernel-coeff} for every $s>0$.

(b) Neither $\mu$ nor $\mu+\nu$ is normalized, and both inequalities in
\eqref{eq:hybrid} are required for $k=0$. In \cite[Lemma~3.1]{Wikstrom2026}
the constant term was kept outside the circle estimate, which forced $\mu$ to
be a probability measure.

(c) For $\nu=0$,~\eqref{eq:hybrid} is the circle-only criterion. For $\mu=0$
it becomes the interpolation problem~\eqref{eq:interpolation}, which
Section~\ref{sec:kernel-only} solves for $c=0.52$ with an absolutely
continuous $\nu$. The circle term, with $\beta_c<1$ inside the annulus,
separates the two sides of~\eqref{eq:hybrid}, so that finitely many atoms in
$\nu$ can work.

(d) Theorem~\ref{thm:hybrid} is the dual form of a linear problem: maximize
$\sum_k(x_k-y_k)/(k+1)$ over non-negative sequences satisfying
\eqref{eq:circle} and~\eqref{eq:kernel-coeff}. Only weak duality is used
below: a sequence with a positive value obstructs the criterion, and
Section~\ref{sec:obstruction} constructs one. We make no claim about the
absence of a duality gap.

(e) Theorem~\ref{thm:kernel} and Proposition~\ref{prop:circle} do not involve
the norm in which $f$ and $g$ are compared. Consider the norm of $L^2$ of a
finite positive radial measure on $\D$ with moments $m_k=\int\abs{z}^{2k}$, so
that $\norm{f}_m^2=\sum_km_k\abs{a_k}^2$. For $f$ and $g$ with finite norm, the
proofs of Proposition~\ref{prop:average} and Theorem~\ref{thm:hybrid} apply
verbatim with $1/(k+1)$ replaced by $m_k$, and so does the obstruction test of
Section~\ref{sec:obstruction}. For the standard weights
$(\alpha+1)(1-\abs{z}^2)^\alpha\,dA$, $\alpha>-1$, the moments are
$m_k=k!\,\Gamma(\alpha+2)/\Gamma(k+\alpha+2)$. Only the explicit interpolant of
Section~\ref{sec:kernel-only} is specific to $1/(k+1)$.
\end{remark}

\section{The certificate for \texorpdfstring{$c_2\ge 0.52755$}{c2 >= 0.52755}}\label{sec:certificate}
\subsection{The data}
Let $c = 0.52755$, $q = c^2$, and $R = 0.99999$ (as exact rational numbers)
and define the measures
\[
 \mu=\sum_{i=1}^{27}w_i\delta_{t_i}+\boldsymbol1_{[R,1]}(t)\,dt,
 \qquad
 \nu=\sum_{j=1}^{749}v_j\delta_{s_j},
\]
where $\delta_t$ denotes a unit point mass at~$t$. The nodes and weights are
terminating decimals, listed in the data file of \cite{Supplement}, and are
treated as exact rational numbers. They satisfy $q\le t_i\le R$, $w_i>0$,
$s_j>0$ and $v_j>0$. The smallest node is $t_1=q$, and the largest is $t_{27} = R$,
with weight $w_R:=w_{27} \approx 4\cdot 10^{-7}$.
The parameters $s_j$ range from $1.49\cdot10^{-782}$ to $0.0427$. For very small
$s$, the weight $h_s(k)$ is close to $1$ for $k$ up to about $\log s/\log q$ and
decays like $q^k/s$ beyond, so these atoms act as mean square comparisons on the
outer circle, truncated at that degree.

Put $\tau_k=R^{k+1}/(k+1)$ and
\[
 J_k=\int_R^1(1-\beta_c(t))\,t^k\,dt.
\]
Since $\int t^k\,d\mu=\sum_iw_it_i^k+(1-R^{k+1})/(k+1)$, the two conditions
\eqref{eq:hybrid} are equivalent to
\begin{align}
 A_k&:=\sum_iw_it_i^k+\hat\nu(k)-\tau_k\ge 0,\label{eq:A}\\
 B_k&:=J_k+\tau_k-\sum_iw_i\beta_c(t_i)t_i^k-\hat\nu(k)\ge 0.\label{eq:B}
\end{align}

\subsection{The terminal interval}
Let
\[
 \kappa_R=\frac\kappa2\,\Phi(qR)\,\Phi(q/R) \qquad\text{and}\qquad
 \gamma_R=\frac{2\kappa_R}{1+\kappa_R(1-R)}.
\]
\begin{lemma}\label{lem:terminal}
For $R\le t<1$ we have $1-\beta_c(t)\ge\gamma_R(1-t)$. Consequently, for
every integer $k\ge 0$,
\begin{equation}\label{eq:J}
 J_k\ge\gamma_R\,\frac{1-R^{k+1}((k+2)-(k+1)R)}{(k+1)(k+2)}.
\end{equation}
\end{lemma}

\begin{proof}
For $p=q^n$ with $n\ge1$, the expression
\[
 \chi(pt)\chi(p/t)=
 \frac{1+p^2-p(t+t^{-1})}{1+p^2+p(t+t^{-1})}
\]
is increasing in~$t$ on $(q,1)$, since $t+t^{-1}$ is decreasing. Hence,
$\Phi(qt)\Phi(q/t)\ge\Phi(qR)\Phi(q/R)$ for $t\ge R$, and since
$1/(1+t)\ge1/2$,~\eqref{eq:beta} gives $\sigma(t)\ge\kappa_R(1-t)$. Since
$1-\beta_c=2\sigma/(1+\sigma)$ is increasing in~$\sigma$,
\[
 1-\beta_c(t)\ge
 \frac{2\kappa_R(1-t)}{1+\kappa_R(1-t)}\ge\gamma_R(1-t).
\]
Integrating $t^k(1-t)$ over $[R,1]$ gives~\eqref{eq:J}.
\end{proof}

\subsection{Enclosures}

For $\Phi_N(x)=\prod_{n=0}^{N-1}\chi(xq^n)$ and $0<x<1$,
\begin{equation}\label{eq:product-tail}
 \Phi_N(x)e^{-\varepsilon_N(x)}\le\Phi(x)\le\Phi_N(x),
\end{equation}
where
\[
 \varepsilon_N(x)=\frac{2xq^N}{(1-xq^N)(1-q)},
\]
since $-\log\chi(v)=2\artanh v\le2v/(1-v)$. We evaluate $\beta_c(t_i)$
through~\eqref{eq:beta} and~\eqref{eq:product-tail} with $N=70$, in Arb ball
arithmetic \cite{Johansson2017}, using \texttt{python-flint}, with~$100$
decimal digits, and use
$\beta_c(q)=1$ for the node $t_1=q$. In $B_k$, the integral $J_k$ is replaced
by its lower bound~\eqref{eq:J}. For $k \ge 1462$ we use the bounds
\[
 0\le\hat\nu(k)\le\Bigl(\sum_jv_j\frac{1+s_j}{s_j}\Bigr)q^k.
\]
For these $k$ the right-hand side is less than $10^{-38}$, so the kernel atoms
are negligible, and in high degrees the certificate rests on the radial part
of~$\mu$ alone. A condition is accepted only if the whole ball is positive. Ball arithmetic
shows $\gamma_R>0.0087522711$.

\subsection{Finitely many degrees}
Table~\ref{tab:certificate} lists the ranges checked in this way, the
smallest certified lower bounds for the relative margins $A_k/\tau_k$ and
$B_k/\tau_k$ (rounded down), and the degrees where they occur. Since $J_k$ is
replaced by~\eqref{eq:J}, and $\hat\nu(k)$ by its bounds for $k\ge 1462$, these
are lower bounds for the true margins.

\begin{table}[ht]
\centering\small
\caption{Finite verification of the certificate at $c=0.52755$.}
\label{tab:certificate}
\begin{tabular}{@{}cccc@{}}
\toprule
Condition & Degrees checked & Smallest certified lower bound & At degree\\
\midrule
$A_k\ge 0$ & $0\le k\le2\,495\,778$ & $2.84806\cdot10^{-8}$ & $k=187\,024$\\
$B_k\ge 0$ & $0\le k\le2\,230\,322$ & $2.84607\cdot10^{-8}$ & $k=7\,235$\\
\bottomrule
\end{tabular}
\end{table}

\subsection{Tails}
\emph{Condition $A_k$.} Let $k_A=\lceil R/w_R\rceil-1=2\,495\,779$. The
exact inequality $w_R(k_A+1)>R$ implies, for $k\ge k_A$,
\[
 \sum_iw_it_i^k\ge w_RR^k\ge\frac{R^{k+1}}{k+1}=\tau_k,
\]
and $\hat\nu(k)\ge 0$, so $A_k\ge 0$.

\emph{Condition $B_k$.}
Let $r_1= 0.999998$ and $r_2=0.999999$ (exact values), so
that $R<r_1<r_2<1$. By Lemma~\ref{lem:terminal},
\[
 J_k\ge\gamma_R\int_{r_1}^{r_2}(1-t)t^k\,dt\ge
 \gamma_R(1-r_2)(r_2-r_1)\,r_1^k.
\]
Moreover,
\[
 \frac{h_s(k+1)}{h_s(k)}=q\,\frac{q^k+s}{q^{k+1}+s}\le R
 \quad\Longleftrightarrow\quad
 q^{k+1}(1-R)\le s(R-q).
\]
Let $k_B=2\,230\,322$. Ball arithmetic shows that
$q^{k_B+1}(1-R)<s_j(R-q)$ for every $j$, and that
\[
 \begin{aligned}
 \gamma_R(1-r_2)(r_2-r_1)\,r_1^{k_B} &> 1.011313\cdot10^{-16},\\
 \sum_iw_it_i^{k_B}+\hat\nu(k_B)     &< 1.011307\cdot10^{-16}.
 \end{aligned}
\]
For $k\ge k_B$, each term of $\sum_iw_it_i^k+\hat\nu(k)$ is multiplied by at
most $R$ when $k$ increases by one, since $t_i\le R$ and by the ratio bound,
whereas $\gamma_R(1-r_2)(r_2-r_1)\,r_1^k$ is multiplied by $r_1>R$. Hence,
$\gamma_R(1-r_2)(r_2-r_1)\,r_1^k>\sum_iw_it_i^k+\hat\nu(k)$ for all
$k\ge k_B$. Since $\beta_c\le1$ and $\tau_k>0$, this gives
$B_k\ge J_k-\sum_iw_it_i^k-\hat\nu(k)>0$. The degree $k_B$ is the first one at
which this comparison holds, which is why the two bounds above are so close.

\begin{proof}[Proof of Theorem~\ref{thm:main}]
Table~\ref{tab:certificate} and the two tail estimates establish
\eqref{eq:A} and~\eqref{eq:B} for every integer $k\ge 0$. Hence, the measures
$\mu$ and $\nu$ satisfy~\eqref{eq:hybrid}, and Theorem~\ref{thm:hybrid}
shows that $\cstar$ is admissible.
\end{proof}

\subsection{How the data were found, and the verification programs}
The data were found by linear programming. The variables were weights on a
grid of radii $t$ and of parameters $\log s$, and the objective was the
smallest relative margin in~\eqref{eq:A} and~\eqref{eq:B}, imposed at all
$k\le1500$ and at $1800$ geometrically spaced degrees up to
$2.8\cdot10^6$. A floating-point scan of every degree up to $3\cdot10^6$
then found a narrow negative trough near $k=1.226\cdot10^6$. Adding
the local minima of the scan as constraints and refitting the weights on the
same support removed it. None of this enters the proof. With Wang's factor
$G_c$ in place of $G'_c$, the same linear programs remain feasible up to
$c\approx0.52754$, but not at $0.52755$, so the sharp factor gains only about
$10^{-5}$ here.

A second program, sharing no code with the first, repeats the verification
with logarithmic product sums and $80$ decimal digits for every
$k\le3\cdot10^6$, closes the tails with the exponential bound for
$\hat\nu(k)$, and obtains the same margins. Data, programs and reports are
archived at~\cite{Supplement}.

\section{The limit of the criterion}\label{sec:obstruction}
By Remark~\ref{rem:criterion}(a), the hypotheses of
Theorem~\ref{thm:hybrid} fail as soon as there is a finitely supported real
sequence $(d_k)$ such that, with $x_k=(-d_k)_+$ and $y_k=(d_k)_+$,
\begin{enumerate}
\item[(i)] $\sum_kd_kh_s(k)\ge 0$ for every $s>0$,
\item[(ii)] $X(t)\le\beta_c(t)Y(t)$ for $q\le t<1$, and
\item[(iii)] $\sum_kd_k/(k+1)<0$.
\end{enumerate}
Indeed, the chain of inequalities in the proof of Theorem~\ref{thm:hybrid}
would then give $\sum_kd_k/(k+1)\ge 0$.

\begin{lemma}\label{lem:monotone}
If a finitely supported sequence satisfies \textup{(i)} and \textup{(ii)}
for some $c\in(0,1)$, then it satisfies \textup{(i)} and \textup{(ii)} for
every $c'\in[c,1)$.
\end{lemma}

\begin{proof}
For (ii), the annulus shrinks as $c$ increases, so $c^*$ and hence $K_c$,
$G'_c$ and $\beta_c$ increase; moreover $[c'^2,1)\subset[c^2,1)$. For (i),
let $q'=c'^2$ and $\alpha=\log q'/\log q\in(0,1]$; we may assume $\alpha<1$.
With $x=q^k$,~\eqref{eq:stieltjes-ratio} represents $q'^k/(q'^k+s)$ as a
positive average of the functions $q^k/(q^k+y)=h_y(k)/(1+y)$, $y>0$, so (i)
at $c$ implies (i) at $c'$.
\end{proof}

\subsection{Reduction to a moment problem}
Fix $N$ and put $\xi_k=q^k$ for $0\le k\le N$. For a sequence supported in
$\{0,\dots,N\}$, partial fractions give
\begin{equation}\label{eq:obstruction-kernel}
 \sum_{k=0}^Nd_k\frac{\xi_k}{s+\xi_k}
 =\frac{P(s)}{\prod_{k=0}^N(s+\xi_k)},
\end{equation}
where $P$ is the real polynomial of degree at most $N$ with
$P(-\xi_k)=d_k\Delta_k$ and $\Delta_k=\xi_k\prod_{l\ne k}(\xi_l-\xi_k)$.
This is a bijection between such sequences and real polynomials $P$ of
degree at most $N$. Thus, (i) holds exactly when $P\ge 0$ on $(0,\infty)$, and
\[
  \sum_{k=0}^N \frac{d_k}{k+1}=\mathcal L(P) \coloneq
  \sum_{k=0}^N \frac{P(-\xi_k)}{(k+1)\Delta_k}.
\]
So the kernel inequalities alone are obstructed at degree $N$ exactly when
the linear functional $\mathcal L$ takes a negative value on a polynomial that
is non-negative on $[0,\infty)$; this is a truncated Stieltjes moment problem
for the moments $\ell_j=\mathcal L(s^j)$. A simple obstruction is
$P(s)=sp(s)^2$. Let $\Lambda_n=(\ell_{i+j+1})_{i,j=0}^{n-1}$. If $\Lambda_9$ is
positive definite, let $p_9$ be the monic polynomial of degree $9$ that is
orthogonal to $1,s,\dots,s^8$ for the form $(p,r)\mapsto\mathcal L(spr)$. Then
\begin{equation}\label{eq:negative-schur}
 \mathcal L(sp_9^2)
 =\ell_{19}-\boldsymbol\ell^{\mathsf T}\Lambda_9^{-1}\boldsymbol\ell
 =\frac{\det\Lambda_{10}}{\det\Lambda_9},\qquad
 \boldsymbol\ell=(\ell_{10},\dots,\ell_{18})^{\mathsf T}.
\end{equation}

\begin{proposition}\label{prop:obstruction}
Let $c_0 = 0.527558$ and $N=20$. Then $\Lambda_9$ is positive
definite and $\mathcal L(sp_9^2)<0$. The sequence $d_k=P(-\xi_k)/\Delta_k$ with
$P(s)=sp_9(s)^2$ satisfies \textup{(i)}, \textup{(ii)} and \textup{(iii)}.
Consequently, the hypotheses of Theorem~\ref{thm:hybrid} fail for every
$c\in[c_0,1)$.
\end{proposition}

\begin{proof}
Throughout the proof $c=c_0$ and $q=c_0^2$. Then $q$, the numbers $\xi_k$
and $\Delta_k$, and the moments $\ell_j$ are rational, and exact elimination
gives nine positive pivots for $\Lambda_9$ and $\mathcal L(sp_9^2)<0$. Thus
(i) holds since $P\ge 0$ on $(0,\infty)$, and (iii) holds since
$\sum_kd_k/(k+1)=\mathcal L(P)$. The coefficients
of $p_9$ are listed in \cite{Supplement}.

For (ii), let $d(t)=\sum_kd_kt^k$. Letting $s\to0$ and $s\to\infty$ in
\eqref{eq:obstruction-kernel} gives $d(1)=P(0)/\prod_k\xi_k=0$ and
$d(q)=\lim_{s\to\infty}sP(s)/\prod_k(s+\xi_k)=0$, since $\deg P=19$. Hence,
$d(t)=(t-q)(1-t)Q(t)$ for a rational polynomial $Q$ of degree $18$. The
first four factors of~\eqref{eq:geo-kappa} give
\[
 \kappa\le\prod_{j=0}^{3}
   \left(\frac{1-cq^j}{1+cq^j}\right)^2<\frac{43}{1000}.
\]
In~\eqref{eq:beta} write $\chi(q/t)=(t-q)/(t+q)$ and bound the other factors
of $\Phi(qt)\Phi(q/t)$ by $1$; with $1-\beta_c=2\sigma/(1+\sigma)\le2\sigma$,
this gives
\[
 1-\beta_c(t)\le\frac{2\kappa(1-t)(t-q)}{(1+t)(t+q)},\qquad q\le t\le1.
\]
The polynomial
\[
 \widetilde Q(t)=(1+t)(t+q)Q(t)-\frac{86}{1000}\,Y(t)
\]
has degree $20$, and all its $21$ coefficients in the Bernstein basis
\[
 \binom{20}{j}\left(\frac{t-q}{1-q}\right)^j
                   \left(\frac{1-t}{1-q}\right)^{20-j},
 \qquad 0\le j\le20,
\]
are positive, by exact computation. Hence, $\widetilde Q\ge 0$ on $[q,1]$, so
$d\ge(1-\beta_c)Y$ there. Since $d=Y-X$, this is (ii).

By Lemma~\ref{lem:monotone}, (i) and (ii) also hold for every
$c\in[c_0,1)$, and (iii) does not depend on $c$.
\end{proof}

\begin{remark}\label{rem:obstruction}
(a) The sequence in Proposition~\ref{prop:obstruction} is constructed from
the kernel inequalities alone; (ii) is verified afterwards. In floating
point, the Hankel test above with $N=30$, $40$ and $60$ fails from
$c\approx0.5275572$ on, close to the value $0.52755$ of
Theorem~\ref{thm:main}. This suggests that the kernel inequalities determine
where the criterion stops. Conditions (i) and (iii) alone show that no
measure $\nu$ satisfies~\eqref{eq:interpolation} when $c\ge 0.527558$, since
integrating (i) against such a $\nu$ would give $\sum_kd_k/(k+1)\ge 0$.

(b) Since Wang's factor satisfies $1-1/G_c\ge\beta_c$, the same sequence
also obstructs the criterion with $G_c$ in place of $G'_c$.

(c) The number $0.527558$ is a limit for the criterion~\eqref{eq:hybrid}.
It is not an upper bound for $c_2$: the sequence $(d_k)$ satisfies the
inequalities used by the criterion, but it is not claimed to come from a
pair $f,g$ satisfying~\eqref{eq:domination}.

(d) By Remark~\ref{rem:positive-kernel}, the hypothesis
$\abs{f}\le\abs{g}$ on $\A$ is equivalent to the positivity of the kernels
$P_s$, and every inequality for the numbers $\abs{a_k}^2$ and $\abs{b_k}^2$
that is obtained by integrating $P_s$ is a positive combination of the kernel
inequalities. The circle estimate adds two-point information, through the
M\"obius pseudodistance, combined with the rank-two bound of
Lemma~\ref{lem:rank-two} on one circle at a time.
Proposition~\ref{prop:obstruction} shows that these two families together
stop at $0.527558$. Any further improvement must therefore use inequalities
that are not consequences of these two families. The correction at the first
Taylor coefficient in~\cite{Pinho2026b}, which strengthens~\eqref{eq:circle}
by a term proportional to $y_1$, is of this kind. Other examples would be
nonlinear consequences of the positivity of $P_s$ that combine several
circles, or inequalities that keep the products $a_j\overline{a_l}$ and
$b_j\overline{b_l}$ with $j\ne l$.
\end{remark}

\section*{Declarations}

\subsection*{Tool and computational resource disclosure}
Claude Opus~5.5 (Anthropic) and GPT-6 Astra (OpenAI) were used as research
assistants for numerical exploration, the search for the certificate, the
verification programs, and drafting. The computations were run on an Apple
M1~Pro computer. The author has checked all arguments and takes full
responsibility for the content.

\subsection*{Data and code availability}
The certificates, the verification programs and their reports are archived
at~\cite{Supplement}.

\end{document}